\documentclass[11pt,reqno]{amsart}
\usepackage{amsmath}
\usepackage{amssymb}
\usepackage{mathrsfs}
\usepackage{mhequ}
\usepackage{color}
\usepackage{graphicx}
\usepackage{enumerate}
\usepackage{todonotes}
\usepackage[unicode=true,
bookmarks=true,bookmarksnumbered=false,bookmarksopen=false,
breaklinks=false,pdfborder={0 0 1},backref=false,colorlinks=true]{hyperref}
\usepackage{bbm} 
\usepackage[small]{caption}
\usepackage{hyphenat}
\usepackage{enumitem}

\newtheorem{maintheorem}{Theorem}

\newtheorem{theorem}{Theorem}[section] 

\newtheorem{definition}[theorem]{Definition}

\newtheorem{remark}[theorem]{Remark}

\newtheorem{example}[theorem]{Example}

\newcommand{\R} {\mathbb R}

\newcommand{\1}{\mathbbm{1}}

\newcommand{\figref}[1]{Figure~\ref{#1}}

\begin{document}
		\title[Characterization of foliations via disintegration maps]%
	{Characterization of foliations via disintegration maps} 
	
	\author[F. M\"unch]{Florentin M\"unch}
	\address{Florentin M\"unch \\
		Department of Mathematics\\
		Universit\"at Leipzig \\
		04109 Leipzig\\
		Germany \\
		\textit{Email address:} \normalfont{\texttt{cfmuench@gmail.com}}
	}

	\author[R. Possobon]{Renata Possobon}
	\address{Renata Possobon\\
		Institute of Mathematics, Department of Applied Mathematics\\
		Universidade Estadual de Campinas
		13.083-859 Campinas - SP\\
		Brazil \\
		\textit{Email address:} \normalfont{\texttt{possobon@ime.unicamp.br}}
	}

	\author[C.S. Rodrigues]{Christian S.~Rodrigues}
	\address{Christian S.~Rodrigues\\
		Institute of Mathematics, Department of Applied Mathematics\\
		Universidade Estadual de Campinas\\
		13.083-859 Campinas - SP\\
		Brazil\\
		and Max-Planck-Institute for Mathematics in the Sciences\\
		Inselstr. 22\\
		04103 Leipzig\\
		Germany \\
		\textit{Email address:} \normalfont{\texttt{rodrigues@ime.unicamp.br}}
	}
	
	\begin{abstract} 
		In this paper, we present a novel approach for analyzing the relationship between the supports of conditional measures and their geometric arrangement in Wasserstein space via the disintegration map. Our method establishes criteria to determine when such conditional measures arise from a metric measure foliation. Additionally, we provide a example demonstrating how this framework can be applied to study perturbations of disintegration-induced foliations.
	\end{abstract}
	
	\maketitle
	
	\section{Introduction}
		
		The disintegration of measures offers a powerful framework for analyzing local and global probabilistic structures by decomposing a measure into conditional probabilities. Roughly speaking, given  measurable spaces $X, Y$, a measure  $\mu$ on $X$, and a measurable map $\pi: X \to Y$, a disintegration allows us to break $\mu$ into a family of conditional measures $\{\mu_y\}$ on $X$, depending measurably on $y$, such that $\mu$ can be reconstructed by integrating $\mu_y$ against $\pi_*\mu$, the push-forward of $\mu$ by $\pi$. This decomposition is crucial in several areas, such as Probability Theory, Ergodic Theory, and Geometric Measure Theory, where it provides a rigorous way to handle conditional distributions. 
		
		Disintegration theorems typically require topological assumptions to ensure existence and uniqueness of such decomposition, making them a key instrument to define measures on fibered spaces or under projections, for instance. Although significant progress has been made in establishing the existence of conditional measures under several different hypothesis, most of the existing approaches neglect the geometric aspects of the underlying space.
		
		In a recent paper \cite{PR24}, the second and third authors used the intrinsic structures of probability spaces to tackle the disintegration of measures from a perspective which takes into account statistical and geometrical properties of dynamical systems. In particular, it introduced the notion of disintegration maps, which are essentially functions that assign to each point $y \in Y$ its corresponding conditional measure $\mu_y$, providing a structured way to analyze a disintegration. We emphasize that these objects can trace geodesics in the space of probability measures, which correspond to the solutions of a time-evolving transport problem, and exhibit certain rigidity phenomena: if one of the endpoints measure in the geodesic path is absolutely continuous, then so are all inner points of the geodesic. Therefore, it becomes evident that certain geometric properties of the base space control the regularity of the conditional measures.
		
		In this paper, we further explore disintegration maps by defining a notion of derivative for them. This derivative resembles a metric derivative, comparing the distance between conditional measures and the distance between the supports of such measures, without requiring a differentiable structure on the base space. Motivated by this derivative, we introduce a novel dense energy functional, similarly to what is done in Geometric Analysis, for cases in which the disintegration applies to every $y \in Y$. This functional characterizes highly ordered geometric structures in the base space: its minimum corresponds precisely to a metric measure foliation, that is, a foliation of the base space whose leaves are parallel to each other and where the distance between conditional measures on these leaves equals the distance between the leaves themselves.
			
		Metric measure foliation arises in various contexts, including Riemannian submersions between weighted manifolds, isometry group actions, and warped products, for instance \cite{GKMS18, Kaz22}. They are significant in Differential Geometry and Analysis due to their role in studying geometric properties, and providing a framework for convergence theory. One the one hand, metric foliations themselves represent a smooth partition of a Riemannian manifold into lower-dimensional submanifolds that are locally equidistant, similar to orbits of an isometric group action or parallel hypersurfaces \cite{GG88}. On the other hand, metric measure foliations are instrumental in the convergence theory of metric measure spaces, for instance. Research indicates that properties like the curvature-dimension condition and the Cheeger energy functional can be preserved from the original space to its quotient via these foliations. This framework also aids in investigating the convergence of sequences of metric measure spaces with unbounded dimensions \cite{Kaz22}.
		
		The main contributions of this paper are presented in Section \ref{main}, where we introduce a notion of derivative for disintegration maps (Definition \ref{def_deriv}) and an energy functional $\mathcal{E}_p$ (Definition \ref{def_energy}), designed to quantify the relationship between the Wasserstein distance of conditional measures and the distance between their supports. Our main result is:
		\begin{theorem}
			Given $f$ a disintegration map, the energy $\mathcal{E}_p (f)$ is equal to one precisely when the conditional measures have supports constituting a metric measure foliation.
		\end{theorem}
		 \noindent The precise statement shall be given in Theorem \ref{e_full}. The paper is organized as follows. In Section \ref{framework}, we establish the theoretical framework and requisite definitions. Section \ref{back} reviews the essential background on disintegration maps. Our core results are presented in Section \ref{main}. These include examples that motivate the hypothesis in Theorem \ref{e_full}, and a further example applying the energy functional to analyze the evolution of a measure under a flow.

	\section{Theoretical framework} \label{framework}
		
		Let $X$, $Y$ be locally compact complete separable metric spaces. Consider a Borel map $\pi: X \to Y$, and $\mu \in \mathcal{M}_{+}(X)$, where $\mathcal{M}_{+}(X)$ is the set of all finite positive Radon measures on $X$. Define $\nu := \pi_{*}\mu$ in $\mathcal{M}_{+}(Y)$. Then, by \cite[Theorem A]{PR24}, there exist measures $\mu_{y} \in \mathcal{M}_{+}(X)$ such that
		
		\begin{enumerate}[leftmargin=*]
			\item $\mu_{y}$ is a probability on $X$ for $\nu$-almost every $y \in Y$;
			\item $\mu_{y}$ is concentrated on $\pi^{-1}(y)$ for $\nu$-almost every $y \in Y$;
			\item $y \mapsto \mu_{y}$ is a Borel map for $\nu$-almost every $y \in Y$; and
			\item $\mu(A)= \int_{Y} \mu_{y} (A) ~d\nu(y)$ for every Borel subset $A$ of $X$.
		\end{enumerate}
		
		The family $\{\mu_y\}$ is called the \textbf{disintegration} of $\mu$ with respect to $\nu$, and the measures $\mu_y$ are called conditional probabilities. In essence, this disintegration theorem provides a specific and constructive way to define conditional probabilities when the conditioning is based on a measurable function. This construction leads to more regular conditional probabilities, a property not universally guaranteed by the more abstract definition of conditional probability with respect to an arbitrary sub-$\sigma$-algebra.

		In order to examine this family of probability measures, we endow the space $\mathscr{P}(X)$ of probability measures on $X$, with a metric structure. The framework of Wasserstein spaces provides a natural geometric setting for such investigations. The \textbf{Wasserstein space} of order $p$, for  $p \in [1, \infty)$, is the set
		\begin{displaymath}
			\mathscr{P}_{p}(X):= \Big\{ \mu \in \mathscr{P}(X) : \int d(x, \tilde{x})^{p} \mu(dx) < +\infty \Big\},
		\end{displaymath}
		with $\tilde{x} \in X$ arbitrary, endowed with the \textbf{Wasserstein distance}
		\begin{displaymath}
			W_{p}(\mu, \nu):= \Big( \inf\limits_{\gamma \in \Pi(\mu, \nu)}  \int d(x_{1}, x_{2})^{p} ~d\gamma (x_{1}, x_{2})\Big)^{\frac{1}{p}},
		\end{displaymath}
		for probability measures $\mu$ and $\nu$ on $X$, where $\Pi(\mu, \nu)$ is the set of measures $\gamma \in \mathscr{P}(X \times X)$  with marginals $\mu$ and $\nu$.  In this context, given a disintegration of $\mu$ with respect to $\nu$, by disintegration map we mean a function that assigns to each $y \in Y$ a conditional measure $\mu_y \in \mathscr{P}_p(X)$.
		
		\
		
		\begin{definition} \cite[Definition 4.1]{PR24} \label{def_disintegration-map}
			Let $X$, $Y$ be locally compact complete separable metric spaces. Consider a Borel map $\pi: X \to Y$, $\mu \in \mathcal{M}_{+}(X)$, and $\nu := \pi_{*}\mu$. Given a disintegration $\{ \mu_y \}$ of $\mu$ with respect to $\nu$, we define the \textbf{disintegration map} by
			\begin{align*}
				f: Y &\to (\mathscr{P}_p(X), W_{p}) \\
				y & \mapsto \mu_{y},
			\end{align*}
			such that $\mu(A)= \int_{Y} f(y) (A) ~d\nu(y)$, for every Borel subset $A$ of $X$. In order to make clear which measures are associated with the disintegration map, we say that ``$f$ is a disintegration map of $\mu$ with respect to $\nu$''.
		\end{definition}

		\
		
		To give an example, consider the unit square $X = [0,1] \times [0,1]$, the two-dimensional Lebesgue measure $\mu = \text{Leb}^2$ on $X$, the projection onto the first coordinate $\pi(x,y) = x$, and the one-dimensional Lebesgue measure $\nu = \text{Leb}^1$ on $[0,1]$. A disintegration of $\mu$ with respect to $\nu$ is a family of measures $\mu_x = \text{Leb}^1 \text{ on } \{x\} \times [0,1]$, with $x \in [0,1]$. Explicitly, for any measurable subset $A \subseteq X$,
		\begin{displaymath}
			\mu(A) = \int_{[0,1]} \mu_x(A) \, d\nu(x) = \int_0^1 \text{Length}\big(A \cap (\{x\} \times [0,1])\big) \, dx,
		\end{displaymath}
		where $\text{Length} (A \cap (\{x\} \times [0,1])) = \text{Leb}^1 (\{ y \in [0,1]  : (x, y) \in A \})$. Moreover, the map $x \mapsto \mu_x(A)$ is measurable,	and $\mu_x$ is supported on $\pi^{-1}(x)$, satisfying the disintegration conditions.
			
		\
		
		\begin{figure}[h!]
			\begin{center}
				\begin{tikzpicture}[scale=1.5]
					\draw (0,0) rectangle (1,1);
					\foreach \x in {0.7} \draw[red, thick] (\x,0) -- (\x,1);
				\end{tikzpicture}
				\caption{Disintegration of $\text{Leb}^2$ with respect to $\text{Leb}^1$: the red line $\{x\} \times [0,1]$ carries $\mu_x = \text{Leb}^1$.}
			\end{center}
		\end{figure}
		
		In this case, the disintegration map is given by
		\begin{align*}
			f: [0,1] &\to (\mathscr{P}_p([0,1] \times [0,1]), W_{p}) \\
			x & \mapsto \mu_{x},
		\end{align*}
		where $\mu_x$ is $\text{Leb}^1$ on $\{x\} \times [0,1]$. This is a basic example that fits into the broader framework of metric foliations, which is the primary focus of our study. Consider a metric space $(X, d)$, and a partition $\mathcal{F}$ of $X$ into closed subsets. We call $\mathcal{F}$ a foliation of $X$, and the elements of this partition are called leaves. If for every $F \in \mathcal{F}$, $F' \in \mathcal{F}$, and every $x \in F$, holds $d(F, F') = d(x, F')$, where $d(F, F') = \inf \{ d(x, x') : x \in F, x' \in F' \}$ and $d(x, F') = d (\{ x \}, F')$, then  $\mathcal{F}$ is called a \textbf{metric foliation}. In case that each leaf is bounded, we say that $\mathcal{F}$ is bounded. Given a metric foliation $\mathcal{F}$ of $X$, define the equivalence relation:
		\begin{equation} \label{rel_quoti}
			x \sim x' ~\iff~ \exists~ F \in \mathcal{F} ~\text{such that}~ x, x' \in F.
		\end{equation}
		Consider the set of equivalence classes under \eqref{rel_quoti} $X^{*} :=X/\sim$ and the projection  $\pi: X \to X^{*}$ onto $X^{*}$. We call $X^{*}$ the quotient space and $\pi$ the quotient map. We define $d^*$ a distance on $X^*$ by 
		\begin{equation}
			d^{*}(y, y') := d(\pi^{-1}(y), \pi^{-1}(y')).
		\end{equation}
		Note that $\pi$ is a submetry: $\pi(B(x, r)) = B(\pi(x), r)$, where $B(x, r)$ is a ball centered at $x$ with radius $r$. Indeed, 
		\begin{align*}
			B(\pi(x), r) =& ~ \{ y \in X^{ *} : d^{*} (y, \pi(x)) < r\} \\
			=& ~\{ y \in X^{ *} : d (\pi^{-1} (y), \pi^{-1}(\pi(x))) < r\} \\
			=& ~\pi(B(x, r)).
		\end{align*}
		
		\noindent When the disintegration of $\mu$ with respect to a metric foliation $\mathcal{F}$ is such that the $p$-Wasserstein distance between two conditional measures $\mu_y$ and $\mu_{y'}$ equals $d^*(y, y')$, we say that $\mathcal{F}$ is a $p$-metric measure foliation. More precisely:
		
		\
		
		\begin{definition} \label{def_mmf}
			Let $X$ be locally compact complete separable metric space, and $\mathcal{F}$ be a metric foliation of $X$. Consider $\mu \in \mathscr{M}_+(X)$, the quotient space $X^*$ as defined above, and the quotient map $\pi: X \to X^*$. We call $\mathcal{F}$ a \textbf{$p$-metric measure foliation}, for $p \in [1, \infty)$, if $\pi_*\mu$ is locally finite Borel measure on $X^*$, and
			\begin{equation}\label{wass_dis_equal}
				W_p(\mu_y, \mu_{y'}) = d^* (y, y')
			\end{equation}
			for any $y, y' \in X^*$, where $\{ \mu_y \}$ is a disintegration of $\mu$ with respect to $\pi_{*}\mu$.
		\end{definition}
		
		\
		
		\begin{remark}
			The general definition for a metric measure foliation requires $p=2$ and \eqref{wass_dis_equal} almost everywhere; see \cite{GKMS18}, for instance. However, our definition is more restrictive, requiring \eqref{wass_dis_equal} everywhere. This restriction is incorporated into the definition in anticipation of the characterization established in Theorem \ref{e_full}.
		\end{remark}
		
		\
		
		Although the definition holds for $p \in [1, \infty)$, our case of interest is $p=2$, since $\mathscr{P}_2(X)$ has a rich geometric structure \cite{GRS25}. Moreover, if $\mathcal{F}$ is a 2-metric measure foliation, then $W_q(\mu_y, \mu_{y'}) = d^* (y, y')$, for any $q \in [1, \infty)$ and every $y, y' \in X^*$, where $\{\mu_y\}$ is the disintegration of $\mu$ with respect to $\pi_*\mu$ \cite[Proposition 3.11]{Kaz22}. Henceforth, the term ``metric measure foliation'' will refer to what we have defined as a ``2-metric measure foliation''.
				
		A very important example of metric measure foliation is related to the action of isometry group \cite{GR23}. Let $(X, d, \mu)$ be a metric measure space and consider a compact topological group $G$. Let
		\begin{displaymath}
			G \times X \ni (g, x) \mapsto gx \in X
		\end{displaymath}
		be an isometric action of $G$ on $X$. Suppose this action is metric measure isomorphic. That is, for every $g \in G$ the map
		$X \ni x \mapsto gx \in X$ is an isometry preserving the measure $\mu$. Consider the $G$-orbit $[x]$ of a point $x \in X$ and the quotient space $X/G$ endowed with the distance
		\begin{displaymath}
			d_{X/G}([x], [x']) = \inf_{g, g' \in G} d(gx, g'x') .
		\end{displaymath}
		Consider the projection map $\pi: X \to X/G$, that is, $\pi$ is given by $x \mapsto [x]$. The family $\mathcal{F} := \{ \pi^{-1}(y) : y \in X/G \}$	is a metric measure foliation on $X$. Other interesting examples arise from Riemannian submersions of weighted Riemannian manifolds; see \cite{GKMS18}, for instance.
		
		\
		
		\
		
		\section{Preliminaries on disintegration maps} \label{back}
		
		In this section, we discuss key properties of disintegration maps, as introduced in Definition \ref{def_disintegration-map} and developed in \cite{PR24}. To start our discussion, suppose $X = Y = \R$, consider  a measurable map $\pi: X \to Y$ and let $\gamma$ be a probability measure supported on the graph of  $\pi$, i.e., $\gamma = (\mathrm{id}, \pi)_* \mu$, where $\mu$ is some probability measure on $X$ and $(\mathrm{id}, \pi)_* \mu$ is the push-forward of $\mu$ under the map $x \mapsto (x, \pi(x))$. Then, the disintegration of $\gamma$ with respect to $\mu$ is given by a family of Dirac measures $\gamma_x = \delta_{\pi(x)}$. Indeed, for every measurable set \(A \subseteq X \times Y\), one has
		\begin{displaymath}
			\gamma(A) = \mu(\{x \in X : (x, \pi(x)) \in A\}).
		\end{displaymath}
		By definition of the disintegration, it holds
		\begin{displaymath}
			\gamma(A) = \int_X \gamma_x(A_x) \, d\mu(x) = \int_X \delta_{\pi(x)}(A_x) \, d\mu(x),
		\end{displaymath}
		with \(A_x = \{y \in Y : (x, y) \in A\}\). Then,
		\begin{displaymath}
			\delta_{\pi(x)}(A_x) = \1_A(x, \pi(x)),
		\end{displaymath}
		and thus:
		\begin{displaymath}
			\gamma(A) = \int_X \1_A(x, \pi(x)) \, d\mu(x),
		\end{displaymath}
		which matches the definition of $\gamma$. This shows that when $\gamma$ is supported on the graph of a function $\pi$, its disintegration is given by $\gamma_x = \delta_{\pi(x)}$. This is a fundamental example in Optimal Transport and Conditional Probability, where deterministic couplings lead to Dirac disintegrations. Observe that when $\pi = \mathrm{id}$, we have a natural isometric immersion of $X$ into $\mathscr{P}_p(X)$, namely the disintegration map $x \mapsto \delta_x$.
		
		In general, we aim to investigate under what conditions geometric information about $X$ can be derived via disintegration, specifically, through the properties of disintegration maps. A key starting point for this analysis is understanding the continuity of these maps. When referring to continuity of the disintegration map, we mean continuity with respect to weak convergence on $\mathscr{P}_p (X)$.	In \cite{PR24} some results are presented regarding the continuity of disintegration maps, namely:
		
		\
			
		\begin{enumerate}[leftmargin=*]
				\item \cite[Propostion 5.2]{PR24} Let $X$, $Y$ be locally compact complete separable metric spaces, and consider $\mu \in \mathcal{M}_{+}(X)$. If a Borel map $\pi: X \to Y$ is such that $\nu := \pi_{*}\mu$ is a Borel measure, then, the disintegration map $f: Y \to (\mathscr{P}_p(X), W_p)$ of $\mu$ with respect to $\nu$ is nearly continuous, for $p \in [1, \infty)$. Here, a map $f$ on an metric space $(Y, \nu)$ is called nearly continuous if for each $\varepsilon > 0$ there exist a closed subset $\mathcal{K} \subset Y$ with $\nu (Y \backslash \mathcal{K}) < \varepsilon$ such that $f$ restricted to $\mathcal{K}$ is continuous.
				
				\
				
				\item \cite[Propostion 5.4]{PR24} Let $X$, $Y$ be locally compact complete separable metric spaces, and consider $\mu \in \mathcal{M}_{+}(X)$. If a Borel map $\pi: X \to Y$ is bijective and $\pi^{-1}$ is continuous, then the disintegration map $f$ of $\mu$ with respect to $\pi_*\mu$ is continuous.
				
				\

				\item \cite[Propostion 5.8]{PR24} Let $X$ be a locally compact complete separable metric space. Consider a metric foliation  $\mathcal{F}$ of $X$, the quotient space $X^*$, and the quotient map $\pi:X \to X^*$. If $\mathcal{F}$ is a metric measure foliation, then the disintegration map of $\mu$ with respect to $\nu := \pi_{*}\mu$ is continuous.
		\end{enumerate}
		
		\
		
		The continuity of disintegration maps enables us to relate paths in the chosen indexing metric space to corresponding paths of measures in Wasserstein space. For instance, in the case of the disintegration map $x \mapsto \delta_x$, there is a direct correspondence between a path in $X$ and the induced path of Dirac measures. In the context of a metric measure foliation, this relationship becomes more nuanced. Here, transversal paths (those crossing the leaves) must satisfy an additional condition concerning the distance between the supports of the conditional measures along the path, a consequence of the parallelism of the leaves.
			
		In the next section we present an approach that allows us to study the relation between the relative position of the conditional measures in the Wasserstein space and the relative position of their respective supports, via disintegration map.
	
		\
		
		\
	
	\section{On foliations and disintegration maps} \label{main}

		In this section, we will assume that the disintegration of measures holds for every point, rather than for $\nu$-almost every point, denoting the family as $\{ \mu_y \}_{y \in Y}$ to explicitly indicate that the conditional measures are defined for every $y \in Y$. In fact, under some structural assumptions, the disintegration can hold for every $y \in Y$. This may occurs, for example, when $\pi$ is a trivial bundle with compact fiber, or when $Y$ is discrete. In such settings, $\mu_y$ can be defined point wise while preserving all disintegration properties.
	
		Motivated by the concept of the metric derivative and the interplay between the relative positions of conditional measures and their supports, we introduce a notion of derivative for disintegration maps.
		
		\
		
		\begin{definition} \label{def_deriv}
			Let $X$ and $Y$ be locally compact complete separable metric spaces. Consider a measure $\mu \in \mathcal{M}_{+}(X)$ and a Borel map $\pi: X \to Y$ such that $\{ \mu_y \}$ is a disintegration of $\mu$ with respect to $\nu = \pi_{*}\mu$. Let $f$ be the respective disintegration map. We write, for $y', y'' \in Y$,
			\begin{displaymath}
				\rho(y', y'') := d(\pi^{-1}(y'), \pi^{-1}(y'')),
			\end{displaymath} 
			and we define a notion of derivative of $f$ by
			\begin{equation} \label{def_gradient}
				| \nabla f(y) |_p := \lim\limits_{\varepsilon \to 0} ~ \sup_{ \substack{\rho(y, y') \leq \varepsilon \\ \rho(y, y'') \leq \varepsilon \\ y' \neq y''}} ~ \frac{W_p (\mu_{y'}, \mu_{y''})}{\rho(y', y'')}
			\end{equation}
			for $p \in [1, \infty)$, where $\mu_{y'} = f(y')$ and $\mu_{y''} = f(y'')$ are conditional measures.
		\end{definition} 
		
		\
		
		We warn the reader that $\rho$ may not satisfy all the distance axioms.
		
		Note that $| \nabla f (y) |_p \geq 1$ everywhere. Indeed, since $\mu_{y'}$ and $\mu_{y''}$ are supported on $\pi^{-1}(y')$ and $\pi^{-1}(y'')$, respectively, we have
		\begin{align*}
			W_{p}^{p}(\mu_{y'}, \mu_{y''}) :=&  \inf\limits_{\gamma \in \Pi(\mu_{y'}, \mu_{y''})}  \int d(x_{1}, x_{2})^{p} ~d\gamma (x_{1}, x_{2}) \\
			\geq & ~ d(\pi^{-1}(y'), \pi^{-1}(y''))^p.
		\end{align*}
		In the case of a $p$-metric measure foliation, the $p$-Wasserstein distance between conditional measures satisfies $W_p (\mu_{y'}, \mu_{y''}) = d(\pi^{-1}(y'), \pi^{-1}(y''))$, which implies $| \nabla f(y) |_p = 1$ for every $y \in Y$. Meanwhile, $|\nabla f(y)|_p = 1$ for every $y \in Y$ does not necessarily imply that the family of fibers $\{\pi^{-1}(y)\}$ is a $p$-metric measure foliation. However, this derivative can be used as a tool to characterize when $\{\pi^{-1}(y)\}$ forms a $p$-metric measure foliation, provided certain additional assumptions are made. Paralleling techniques from Geometric Analysis and Optimal Transport, we introduce a dense energy functional designed to further analyze the relationship between the distances of conditional measures and their supports. 
		
		\
		
		\begin{definition} \label{def_energy}
			Let $X$ and  $Y$ be locally compact complete separable metric spaces. Consider a measure $\mu \in \mathcal{M}_{+}(X)$ and a Borel map $\pi: X \to Y$ such that $\{ \mu_y \}_{y \in Y}$ is a disintegration of $\mu$ with respect to $\nu = \pi_{*}\mu$. Let $f$ be the respective disintegration map. The $p$-energy of $f$ is given by
			\begin{displaymath}
				\mathcal{E}_p(f) := \| \nabla f \|_{\infty, p} =  \sup_{y \in Y} | \nabla f(y) |_p.
			\end{displaymath}
		\end{definition}
		
		\

		This energy definition involves analyzing $ | \nabla f(y) |_p$ at every point $y \in Y$. Example \ref{everywhere} demonstrates a scenario in which the definition is satisfied for almost every point, yet the characterization we seek to prove does not hold. In the following theorem, we relate the value of this energy functional to settings in which the disintegration corresponds to a metric measure foliation case.

		\begin{maintheorem} \label{e_full}
			Let $(X, d_X)$ be a geodesic locally compact complete separable metric space. Consider a locally compact complete separable metric space $Y$, a measure $\mu \in \mathcal{M}_+(X)$ and a Borel map $\pi: X \to Y$ such that $\{ \mu_y \}_{y \in Y}$ is a disintegration of $\mu$ with respect to $\nu := \pi_{*}\mu$, with $\operatorname{supp}(\mu_y) = \pi^{-1}(y)$. Let $f$ be the respective disintegration map. Then, $\mathcal{E}_p (f) = 1$ if and only if $\{\pi^{-1}(y)\}$ defines a $p$-metric measure foliation on $X$.
		\end{maintheorem}
		
		\
	
		\begin{remark}
			The assumptions on the space $Y$ are made only to guarantee the existence of the disintegration of $\mu$ with respect to $\pi_*\mu$, according \cite[Theorem A]{PR24}. Indeed, the claim of our theorem only concerns the distance between the fibers, and not the original distance on $Y$.
		\end{remark}
		
		\
		
		\begin{proof}
			The demonstration will be carried out in three steps. Step 1 shows that $p$-metric measure foliation implies $\mathcal{E}_p (f) = 1$, and Step 2 and Step 3 show that if $\mathcal{E}_p (f) = 1$, then $\{\pi^{-1}(y)\}$ defines a $p$-metric measure foliation on $X$. More specifically, Step 2 shows that the Wasserstein distance between two conditional measures coincides with the distance between the corresponding fibers; Step 3 shows the metric foliation property, that is, all points on a fiber have the same distance to another fiber.
				
				\
				
			\noindent In the sequel, we consider the distance between the preimages of $y'$ and $y''$ to be given by 
			$\rho(y', y'') = d(\pi^{-1}(y'), \pi^{-1}(y'')) = \inf \{ d_X (x', x'') : \mbox{$x' \in \pi^{-1}(y')$}$, $x'' \in \pi^{-1}(y'') \}$.
				
				\
				
				\
				
				\noindent \textbf{Step 1.} Suppose that $\{\pi^{-1}(y)\}$ defines a $p$-metric measure foliation on $X$. Then, for any $y', y'' \in Y$, holds $W_p(\mu_{y'}, \mu_{y''}) = d (\pi^{-1}(y'), \pi^{-1}(y''))$. Therefore, for any $y \in Y$, one has
				\begin{displaymath}
					| \nabla f(y) |_p = \lim\limits_{\varepsilon \to 0} ~ \sup_{ \substack{\rho(y, y') \leq \varepsilon \\ \rho(y, y'') \leq \varepsilon \\ y' \neq y''}} ~ \frac{W_p (\mu_{y'}, \mu_{y''})}{\rho(y', y'')} = 1,
				\end{displaymath}
				which implies $\mathcal{E}_p (f) = 1$.
				
				\
				
				\

				\noindent \textbf{Step 2.} Suppose $\mathcal{E}_p (f) = 1$. We want to show that in this case, for any $y', y'' \in Y$, it holds that $W_p(\mu_{y'}, \mu_{y''}) = d (\pi^{-1}(y'), \pi^{-1}(y''))$.

				\
				
				\noindent Consider $\varepsilon> 0$ and $y_0, y_1 \in Y$. Then there exist $x_0 \in \pi^{-1}(y_0)$, $x_1 \in \pi^{-1}(y_1)$, such that $d_X(x_0, x_1) \leq d(\pi^{-1}(y_0), \pi^{-1}(y_1)) + \varepsilon$. Let $\gamma: [0, 1] \to X$ be a constant speed minimizing geodesic such that $\gamma(0) = x_0$ and $\gamma(1) = x_1$. Define $y_t := \pi (\gamma(t))$, $t \in [0, 1]$. Note that, for every $s, t \in [0, 1]$, it holds 
				\begin{equation} \label{Wineq_1}
					W_p (\mu_{y_s}, \mu_{y_t}) \geq d(\pi^{-1}(y_s), \pi^{-1}(y_t)) 
				\end{equation}
				and
				\begin{equation} \label{Wineq_2}
					d(\pi^{-1}(y_s), \pi^{-1}(y_t)) \leq |s-t| d_X(x_0, x_1).
				\end{equation}
				
				\
				
				\noindent Define $\varphi(t):= W_p ( \mu_{y_t}, \mu_{y_0})$. We will show that $\varphi$ is continuous and that
				\begin{displaymath}
					\partial_t^+ (\varphi (t)) := \limsup_{h \to 0} \frac{\varphi (t+h) - \varphi (t) }{h} \leq d_X (x_0, x_1).
				\end{displaymath}
				Indeed,
					\begin{align*}
						\limsup_{h \to 0} \left| \frac{\varphi (t+h) - \varphi (t) }{h} \right| =& \limsup_{h \to 0} \left| \frac{W_p (\mu_{y_{t+h}}, \mu_{y_0}) - W_p (\mu_{y_{t}}, \mu_{y_0})}{h} \right| \\
						\leq & \limsup_{h \to 0} \frac{W_p (\mu_{y_{t+h}}, \mu_{y_t})}{|h|} \\
						\leq & d_X(x_0, x_1) \limsup_{h \to 0} \frac{W_p (\mu_{y_{t+h}}, \mu_{y_t})}{d (\pi^{-1}(y_{t+h}), \pi^{-1}(y_t))} \\
						\leq & d_X(x_0, x_1) |\nabla f (y_t)|_p   \\
						\leq & d_X(x_0, x_1) \mathcal{E}_p(f) \\
						=& d_X(x_0, x_1). 
					\end{align*}				
				Here, the first inequality follows from the triangle inequality for the Wasserstein distance; the second one follows from \eqref{Wineq_2}; the third one follows from the definition of the gradient \eqref{def_gradient}; the last equation follows from the assumption $\mathcal{E}_p(f) = 1$. Now, the continuity of $\varphi$ follows easily.
				
				\
				
				\noindent By \cite[Theorem 2]{HT06}, we can integrate the inequality $\partial_t^+ (\varphi (t)) \leq d_X (x_0, x_1)$, giving $\varphi(1) - \varphi(0) \leq d_X (x_0, x_1)$. Thus,
				\begin{displaymath}
					W_p (\mu_{y_0}, \mu_{y_1}) \leq d_X (x_0, x_1) \leq  d(\pi^{-1}(y_0), \pi^{-1}(y_1)) + \varepsilon.
				\end{displaymath}
				Since $\varepsilon$ is arbitrary, taking the limit $\varepsilon \to 0$, we have
				\begin{equation} \label{ineq_2}
					W_p (\mu_{y_0}, \mu_{y_1}) \leq  d(\pi^{-1}(y_0), \pi^{-1}(y_1)).
				\end{equation}
				
				\
				
				\noindent By \eqref{Wineq_1} and \eqref{ineq_2} we obtain $W_p(\mu_{y'}, \mu_{y''}) = d (\pi^{-1}(y'), \pi^{-1}(y''))$.

				\
				
				\
				
				\noindent {\textbf{Step 3.}} Keep assuming $\mathcal{E}_p (f) = 1$. We will show that for every $y_0, y_1 \in Y$, and for every $x_0 \in \pi^{-1}(y_0)$, it holds $d(x_0, \pi^{-1}(y_1)) = d(\pi^{-1}(y_0), \pi^{-1}(y_1))$. 
				
				\
				
				\noindent Note that, by definition, one has  $d(x_0, \pi^{-1}(y_1)) \geq d(\pi^{-1}(y_0), \pi^{-1}(y_1))$. 
				
				\
				
				\noindent Consider a function $r: \pi^{-1}(y_0) \to \R^+_0$ given by $r(x) := d(x, \pi^{-1}(y_1))$. Note that $r(x) \geq d(\pi^{-1}(y_0), \pi^{-1}(y_1))$ for every $x \in \pi^{-1}(y_0)$ and $r$ is continuous. Moreover, note that for $(x, x') \in \pi^{-1}(y_0) \times \pi^{-1}(y_1)$, one has $d_X (x, x') \geq r(x)$. Let $\eta$ be an optimal transport plan from $\mu_{y_0}$ to $\mu_{y_1}$. Then, by Step 2, it holds
				\begin{displaymath}
					 \int d_X (x, x')^p d \eta = W^p_p (\mu_{y_0}, \mu_{y_1}) =  [d(\pi^{-1}(y_0), \pi^{-1}(y_1))]^p.
				\end{displaymath}
				Furthermore, one has
				\begin{displaymath}
					\int d_X (x, x')^p d \eta \geq \int r(x)^p d \mu_{y_0} \geq [d(\pi^{-1}(y_0), \pi^{-1}(y_1))]^p.
				\end{displaymath}
				Since $r^p$ is continuous and $\operatorname{supp}(\mu_{y_0}) = \pi^{-1}(y_0)$, for every $x \in \pi^{-1}(y_0)$ we have $r(x)^p = [d(\pi^{-1}(y_0), \pi^{-1}(y_1))]^p$ . Therefore, for every $x_0 \in \pi^{-1}(y_0)$, we have $d(x_0, \pi^{-1}(y_1)) = d(\pi^{-1}(y_0), \pi^{-1}(y_1))$ 
				
				\
				
				\noindent Step 2 and Step 3 complete the proof.

		\end{proof}
		
		\
		
		In Theorem \ref{e_full}, we characterize a foliation on $X$ via the energy of a disintegration map. We emphasize that the hypotheses of conditional measures having full support and the analysis of the energy functional over all points are fundamental to the result. In the following examples, we present cases that justify the use of such assumptions. 
			
		\
		
		\begin{example} \label{everywhere}
		We will show that there exists a disintegration $\{ \mu_y \}$, such that $|\nabla f (y)|_p = 1$ for $\nu$-almost every $y$, but the metric measure foliation condition does not hold. Consider $X = [-2, 1]$, $Y=[0, 1]$ and $\pi: X \to Y$ such that
		\begin{displaymath} 
			\pi(x) =
			\begin{cases}
				x &, \text{if} ~ x \geq 0 \\
				(\phi + \text{id})^{-1} (-x) &, \text{if} ~ x < 0
			\end{cases}
		\end{displaymath}
		where $\phi$ is the Cantor ternary function. Consider a measure $\mu$ on $X$ such that $\nu = \pi_* \mu$ is a uniform measure on $Y$ and there exists a disintegration of $\mu$ with respect to $\nu$ such that the conditional measures are
		\begin{displaymath}
			\mu_y = \frac{1}{2} \delta_{y} + \frac{1}{2} \delta_{\varphi (y)}
		\end{displaymath}
		where $\varphi (y) = \phi (y) + y$. Note that, for $y', y'' \in Y$, 
		\begin{align*}
			d(\pi^{-1}(y'), \pi^{-1}(y'')) =& d (\{ y', \varphi(y') \}, \{ y'', \varphi(y'') \}) \\
			=& \min \left\{ |y' - y''|, |\varphi(y'), \varphi(y')| \right\} \\
			=& |y' - y''|.
		\end{align*}
		Moreover, let $\mathcal{U}$ be the union of the maximum intervals where the Cantor function is constant, i.e., $\mathcal{U} = (\frac{1}{3}, \frac{2}{3}) \cup (\frac{1}{9}, \frac{2}{9}) \cup (\frac{7}{9}, \frac{8}{9}) \cup \dots$. Note that $\mathcal{U}$ is dense in $[0, 1]$, $\nu (\mathcal{U}) = 1$, and for $y \in \mathcal{U}$, for all $y', y'' \in B_{\varepsilon}(y)$, $\varepsilon < d_Y (y, \mathcal{U}^c)$, we have $W_p(\mu_{y'}, \mu_{y''}) = |y' - y''|$. Therefore, $| \nabla f(y) |_p = 1$ for $\nu$-almost every $y \in Y$.
				
		\end{example}
		
		\
		
		Example \ref{everywhere} justifies the hypothesis of taking the supremum over all $y \in Y$ in Definition \ref{def_energy}. If we had defined the energy functional using an essential supremum, for example, pathological cases like this would arise, where the other hypotheses of Theorem \ref{e_full} are satisfied, but we do not have a $p$-metric measure foliation. Another important remark: when considering a general disintegration $\{ \mu_y \}$, which is defined $\nu$-almost everywhere, we could define measures at the missing points $y$, but this does not alter the geometric structure on $X$. Moreover, note that by requiring that the conditional measures have full support, the fiber structure becomes detectable by the energy functional. The following example illustrates a case in which Theorem \ref{e_full} does not hold if the full support condition is not assumed.

		\
		
\begin{figure}[h!]
			\centering \includegraphics[scale=0.4]{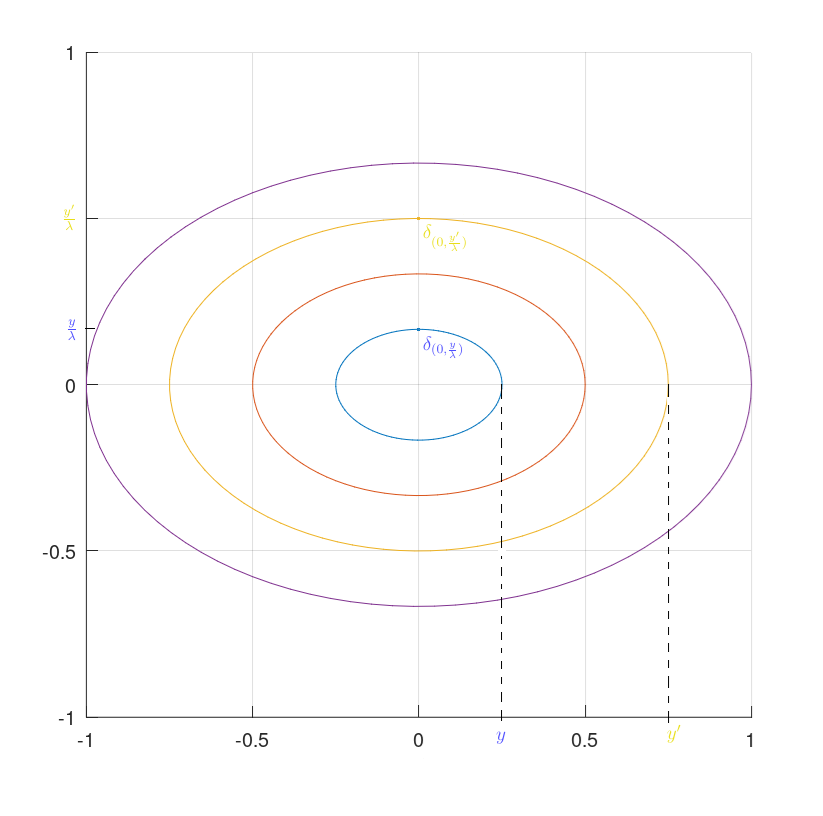}
			\caption{Illustration for the case in which the conditional measures do not have full support: $d(\pi^{-1}(y), \pi^{-1}(y')) = W_p(\mu_{y'}, \mu_{y''})$ for every $y', y'' \in Y$, but $\{ \pi^{-1}(y) \}$ (the ellipses) are not a metric foliation.} \label{elipses_delta}
		\end{figure}

		\begin{example}
			Consider $X = [-1,1] \times \mathbb{R}$, $Y=[0, 1]$, $\pi: X \to Y$ such that $\pi^{-1}(y)$ is an ellipse given by $x_1^2 + \lambda^2 x_2^2 = y^2$, $\lambda >1$, with major radius $y$. Consider a measure $\mu \in \mathcal{M}_{+}(X)$ such that the disintegration of $\mu$ with respect to $\nu$ is $\mu_y = \delta_{(0, \frac{y}{\lambda})}$. In this case,
			\begin{displaymath}
				d(\pi^{-1}(y), \pi^{-1}(y')) = \left| \frac{y}{\lambda} - \frac{y'}{\lambda} \right| = W_p(\mu_y, \mu_{y'})
			\end{displaymath}
			for every $y, y' \in Y$, but $\{ \pi^{-1}(y) \}$ does not define a metric foliation. A representation of the case $\lambda = 1.5.$ is provided in Figure \ref{elipses_delta}.
		\end{example}
		
		\
		
		In the next example, we analyze the disintegration of a measure $\mu$ under a flow $\varphi$. For $t=0$, we have the characterization of a metric measure foliation given by Theorem \ref{e_full}. We use the energy functional to analyze the change in the foliation associated with the measure $\mu_t$, which is the push-forward of $\mu$ by the flow $\varphi$ at time $t$. This example shows that the energy functional exhibits a certain sensitivity to smooth changes in the foliation structure.

		\
		
		\begin{example} \label{ex_eli}
			Consider $X=\{ (x_1, x_2) \in \mathbb{R}^2 ~:~ x_1^2 + x_2^2 \leq R^2 \}$, $R>0$, $\mu \in \mathscr{P}(X)$ with full support and uniformly distributed on $X$, and $Y =(0, R]$. Consider a foliation of $X$ into circles, that is, consider a foliation $\mathcal{F}$ of $X$ such that $F_y = \{ (x_1, x_2) \in \mathbb{R}^2 ~:~ x_1^2 + x_2^2 = y^2 \}$, for each $y \in Y$, and $\pi: X \to Y$ such that $\pi ((x_1, x_2)) = \sqrt{x_1^2 + x_2^2} = y$. Let $\tilde{f}$ be the disintegration map of $\mu$ with respect to $\pi_*\mu$. We will show that $\mathcal{E}_1 (\tilde{f}) = 1$.
			
			\
			
			\noindent Consider this setup in polar coordinates. Since $\mu$ is a uniformly distributed probability measure on $X$, we can write:
			\begin{displaymath}
				d\mu = \frac{dA}{A} = \frac{r ~dr ~d\theta}{\pi R^2}.
			\end{displaymath}
			Define $\nu := \pi_{*} \mu$. Note that, given a open interval $(a, b) \subset (0, R]$, $\nu ((a, b))$ is the area between circles of radii $a$ and $b$, normalized by $A$:
			\begin{displaymath}
				\nu ((a, b)) = \frac{\pi b^2 - \pi a^2}{\pi R^2} = \frac{b^2 - a^2}{R^2}.
			\end{displaymath}
			In this way,
			\begin{displaymath}
				d\nu (r) = \frac{2r}{R^2} ~dr.
			\end{displaymath}
			Let $\{ \mu_y \}_{y \in Y}$ be the disintegration of $\mu$ with respect to $\nu$, such that each $\mu_y$ is a uniformly distributed probability measure on $F_y$. Due to the symmetry of the circle, the arc length is proportional to $\theta$, and 
			\begin{displaymath}
				d\mu_y (r, \theta) = \frac{1}{2\pi} ~ \delta_y ~ d\theta. 
			\end{displaymath}
			Note that
			\begin{displaymath}
				\int_Y \mu_y ~d\nu = \int_{0}^{R} \int_{0}^{2 \pi} \Big( \frac{1}{2 \pi} \delta_y d\theta \Big) \frac{2r}{R^2} ~dr = \frac{1}{\pi R^2} \iint r ~dr ~d\theta,
			\end{displaymath}
			so we recover the measure $\mu$ by integrating the family $\mu_y$ with respect to $\nu$. Moreover, $d(\pi^{-1}(y), \pi^{-1}(y')) = |y - y'|$. One way to write the 1-Wasserstein distance is
			\begin{displaymath}
				W_1 (\mu_y, \mu_{y'}) = \inf_{T} \int d(x, T(x)) ~d\mu_y
			\end{displaymath}
			and since we have measures uniformly distributed and symmetry in the leaves, the function $T$ that minimizes this problem is such that, for every $x \in F_y$, $T(x)$ is a point in $F_{y'}$ such that $d(x, T(x)) = |y - y'|$, that is, the transport occurs radially, orthogonally to the leaves. Therefore,
			\begin{displaymath}
				W_1 (\mu_y, \mu_{y'}) = \int_{F_y} |y - y'| ~d\mu_y = |y - y'|.
			\end{displaymath}
			Then, $\mathcal{E}_1 (\tilde{f}) = 1$.
			
			\
			
			In this case, we have a $1$-metric measure foliation. Suppose that $\mu$ evolves under a flow $\varphi$. If this flow spirals toward the origin, for example, the circular leaves are preserved under the action of the flow, with the radius shrinking. The disintegration structure is preserved by the flow, and the energy remains constant for each disintegration map $\tilde{f_t}$ associated with the disintegration of $\mu_t = (\varphi_t)_*\mu$ with respect to $\nu$. However, if the flow implements, for example, a vertical compression, the circular leaves become ellipses with increasing eccentricity over time. In this case, the disintegration structure changes in time. We will show that the value of $\mathcal{E}_p(\tilde{f_t})$  varies proportionally with this change in the disintegration structure.
			
			\			
			
			\noindent Consider $X=\{ (x_1, x_2) \in \mathbb{R}^2 ~:~ x_1^2 + \lambda^2 x_2^2 \leq R^2 \}$, $R>0$, $\lambda \geq 1$, $\mu \in \mathscr{P}(X)$ with full support and uniformly distributed on $X$, $Y =(0, R]$, and a foliation $\mathcal{F}$ of $X$ such that $F_y = \{ (x_1, x_2) \in \mathbb{R}^2 ~:~ x_1^2 + \lambda^2 x_2^2 = y^2 \}$, for each $y \in Y$, and $\pi: X \to Y$ such that $\pi ((x_1, x_2)) = \sqrt{x_1^2 + \lambda^2 x_2^2} = y$. In polar coordinates
			\begin{displaymath}
				d\mu = \frac{dA}{A} = \frac{\frac{r}{\lambda} ~dr ~d\theta}{ \frac{\pi R^2}{\lambda} } = \frac{r ~dr ~d\theta}{\pi R^2}.
			\end{displaymath}
			Define $\nu := \pi_{*} \mu$. Note that, given a open interval $(a, b) \subset (0, R]$, $\nu ((a, b))$ is the area between ellipses of major axis $a$ and $b$, normalized by $A$:
			\begin{displaymath}
				\nu ((a, b)) = \frac{\frac{\pi b^2}{\lambda} - \frac{\pi a^2}{\lambda}}{\frac{\pi R^2}{\lambda}} = \frac{b^2 - a^2}{R^2}.
			\end{displaymath}
			In this way,
			\begin{displaymath}
				d\nu (r) = \frac{2r}{R^2} ~dr.
			\end{displaymath}
			Let $\{ \mu_y \}_{y \in Y}$ be the disintegration of $\mu$ with respect to $\nu$, such that each $\mu_y$ is a uniformly distributed probability measure on $F_y$. However, note that in this case the relation arc length and $\theta$ is no longer linear (see \figref{normalized_arc}).

			\begin{figure}[h!]
				\centering \includegraphics[scale=0.4]{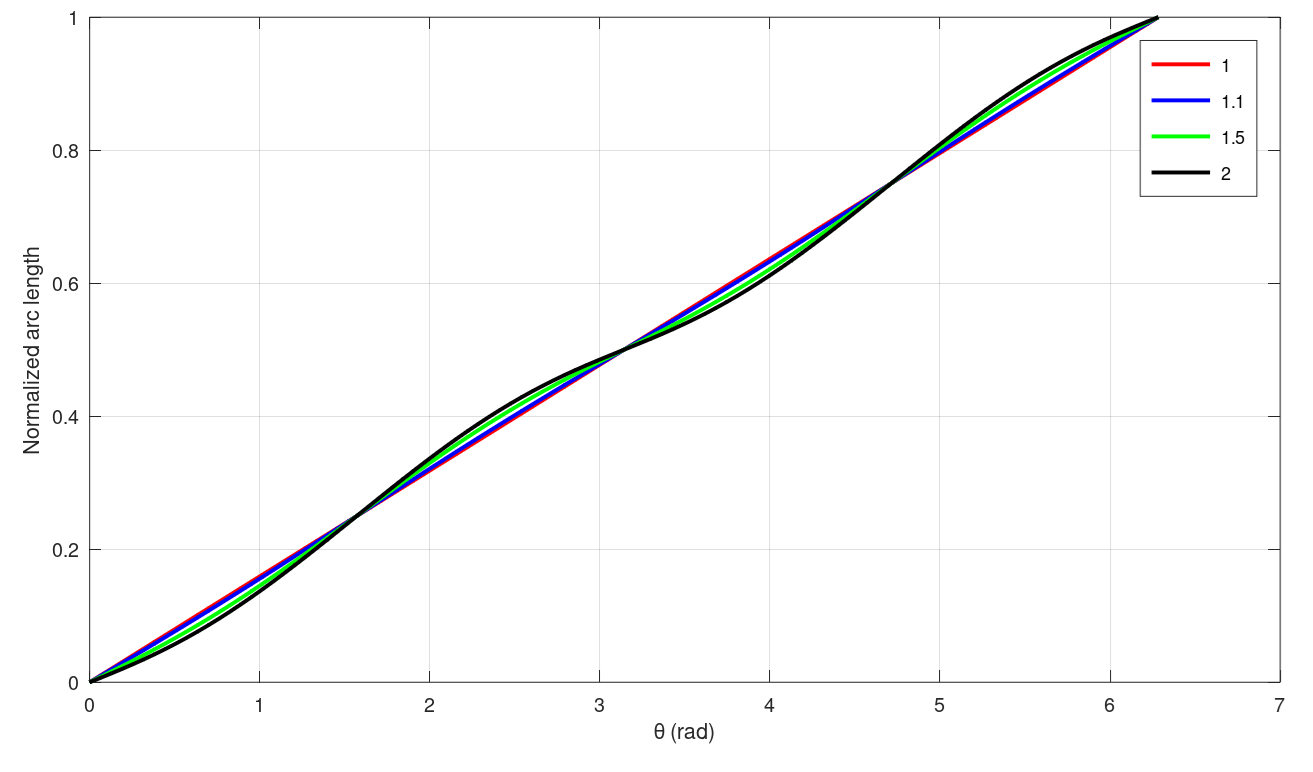}
				\caption{Normalized arc length in function of $\theta$ for different $\lambda$.} \label{normalized_arc}
			\end{figure}
			
			\

			\noindent We have
			\begin{displaymath}
				d\mu_y (r, \theta) = \frac{y ~ \sqrt{1 - (1 - \lambda^{-2}) \cos^2 (\theta)}}{\mathcal{L}_y} ~\delta_y ~d\theta.
			\end{displaymath}
			In this case, $d(\pi^{-1}(y), \pi^{-1}(y')) = \frac{1}{\lambda}|y - y'|$ and, to calculate the $1$-Wasserstein distance between $\mu_y$ and $\mu_{y'}$, we still have the radial transport, although $d(x, T(x))$ depends of $\theta$: for $x = y \cos (\theta)$, $d(x, T(x)) =  R_{y'} (\theta) - R_y (\theta)$, where
			\begin{displaymath}
				R_{y} (\theta) = \frac{y}{\sqrt{\lambda^2 - (\lambda^{2} - 1) \cos^2 (\theta)}}
			\end{displaymath}
			and		 
			\begin{displaymath}
				R_{y'} (\theta) = \frac{y'}{\sqrt{\lambda^2 - (\lambda^{2} - 1) \cos^2 (\theta)}}.
			\end{displaymath}
			Then,
			\begin{align*}
				W_1 (\mu_y, \mu_{y'}) =& \int  R_{y'} (\theta) - R_y (\theta) ~d\mu_y \\
				=& \int \frac{y' - y}{\sqrt{\lambda^2 - (\lambda^{2} - 1) \cos^2 (\theta)}} \frac{y ~ \sqrt{1 - (1 - \lambda^{-2}) \cos^2 (\theta)}}{\mathcal{L}_y} ~\delta_y ~d\theta \\
				=& \frac{(y'-y)~y}{\lambda ~\mathcal{L}_y} ~ \int_{0}^{2\pi} ~d\theta \\
				=& \frac{2 \pi ~(y'-y)~y}{\lambda ~\mathcal{L}_y}.
			\end{align*}
			
			\
			
			\noindent Let $\tilde{f}_{\lambda(t)}$ be the disintegration map of $\mu_t$ with respect to $\nu$. Based on the developments above, we have that

			\begin{displaymath}
				\mathcal{E}_1(\tilde{f}_{\lambda(t)}) = \frac{2 \pi y}{\mathcal{L}_y},
			\end{displaymath}
			that is, the energy is the ratio between the length of the circle with radius $y$ and the length of the ellipse with major axis equal to $y$. Thus, we have that as $\mu$ evolves over time, the structure for the disintegration of $\mu_t$ with respect to $\pi_*\mu_t$ also changes, and this variation is captured by the energy functional. Indeed, Figure \ref{e_1_lambda} depicts the dependence of $\mathcal{E}_1(\tilde{f}_{\lambda(t)})$ on $\lambda(t)$. As shown, $\mathcal{E}_1$ grows smoothly with increasing foliation perturbation.

			\
			
			\begin{figure}[h!]
				\centering \includegraphics[scale=0.38]{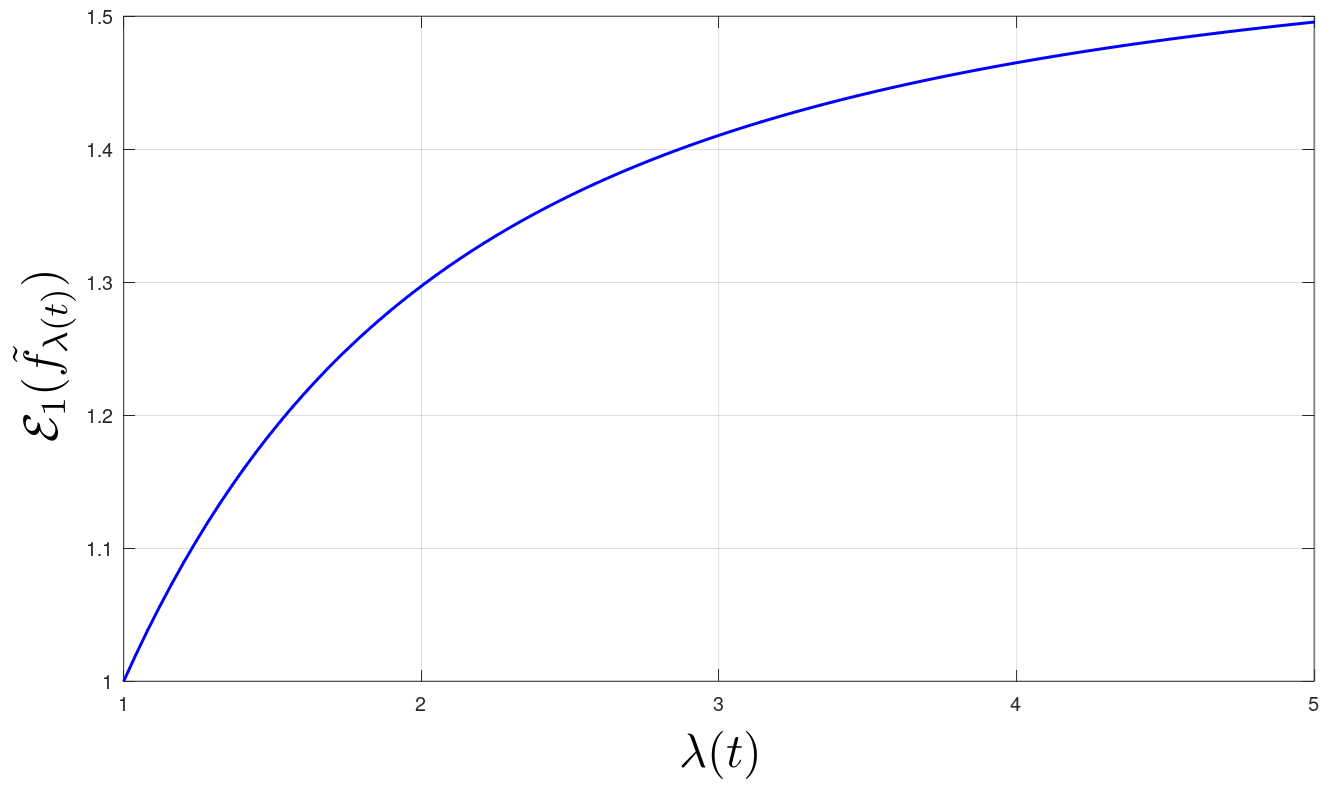}
				\caption{Energy as a function of the eccentricity of the leaves.} \label{e_1_lambda}
			\end{figure}
			
			\

			\noindent For comparison purposes, Table \ref{comp} displays the values of $\mathcal{L}_y$ for $\lambda = 1$ (circle), $\lambda = 1,001$, $\lambda = 1,01$, $\lambda = 1,1$, $\lambda = 1,5$ and $\lambda = 2$. Table \ref{table_energy} summarizes the computed values of $\mathcal{E}_1(\hat{f}_{\lambda})$ corresponding to the cases in Table \ref{comp}. 
			\
			
			\
			
			\begin{table}[h!]
				\caption{$\mathcal{L}_y$ for different values of $\lambda$} \label{comp}
				\footnotesize
				\begin{tabular}{ c | c c c c c c}
					y & $\lambda = 1$ & $\lambda = 1,001$ & $\lambda = 1,01$ & $\lambda = 1,1$ & $\lambda = 1,5$ & $\lambda = 2$ \\
					\hline
					$1$ & $6,28318531$ & $6,28004725$	& $6,25211912$ & $6,00098645$ & $5,28847986$ &	$4,84422411$
				\end{tabular}
			\end{table}
			
			\
			
			\begin{center}
				\begin{table}[h!]
					\caption{Values of $\mathcal{E}_1(\hat{f}_{\lambda})$} \label{table_energy}
					\begin{tabular}{c c c c c c}
						$\lambda = 1$ & $\lambda = 1,001$ & $\lambda = 1,01$ & $\lambda = 1,1$ & $\lambda = 1,5$ & $\lambda = 2$ \\
						\hline
						$1$ & $1,000499687$	& $1,004968906$ & $1,047025412$	& $1,188089106$ & $1,297046785$
					\end{tabular}
				\end{table}
			\end{center}
		\end{example}

		\
		
		Example \ref{ex_eli} highlights the versatility of our approach: beyond characterize foliations on $X$ via disintegration maps, it also provides a framework for analyzing structural perturbations in foliations. The concepts developed here have potential applications across Dynamical Systems, Machine Learning, and Stochastic Analysis.

	\section*{Acknowledgements}
	\noindent This work was financed in part by the Coordena\c{c}\~ao de Aperfei\c{c}oamento de Pessoal de N\'ivel Superior - Brasil (CAPES) - Finance Code 001. C. S. R. has been partially supported by S\~{a}o Paulo Research Foundation (FAPESP): grant \#2016/00332-1, grant \#2018/13481-0, and grant \#2020/04426-6. The opinions, hypotheses and conclusions or recommendations expressed in this work are the responsibility of the authors and do not necessarily reflect the views of FAPESP.
	
	\noindent C. S. R. would like to acknowledge support from the Max Planck Society, Germany, through the award of a Max Planck Partner Group for Geometry and Probability in Dynamical Systems. R. P. would like to acknowledge support from the Max Planck Institute for Mathematics in the Sciences and from CAPES.
	
	\noindent Our manuscript has no associated data.
	
	\noindent The authors declare no conflict of interest.

	\bibliographystyle{amsalpha}
	
\end{document}